\documentclass[oneside, reqno] {amsart}

\usepackage{xypic}
\input xy
\xyoption{all}
\usepackage{epsfig}
\usepackage{amsthm}
\usepackage{amssymb}
\usepackage{amsmath}
\usepackage{amscd}
\usepackage{color}
\usepackage[T1]{fontenc}
\usepackage{upgreek}
\usepackage[font=scriptsize]{caption}
\usepackage{wrapfig}

\usepackage{graphicx}
\usepackage{subfigure}

\newcommand{\bg}{\begin{equation}}
\newcommand{\ed}{\end{equation}}
\newcommand{\bga}{\begin{eqnarray}}
\newcommand{\eda}{\end{eqnarray}}
\newcommand{\pf}{\textbf{Proof:\ }}

\def\cbdu{\par{\raggedleft$\Box$\par}}

\newtheorem {Theorem}  {Theorem}

\numberwithin{Theorem}{section}

\newtheorem {Lemma}[Theorem]  {Lemma}

\theoremstyle{definition}

\theoremstyle{remark}
\newtheorem{Remark}[Theorem]{\bf Remark}

\expandafter\chardef\csname pre amssym.def
at\endcsname=\the\catcode`\@ \catcode`\@=11
\def\undefine#1{\let#1\undefined}
\def\newsymbol#1#2#3#4#5{\let\next@\relax
 \ifnum#2=\@ne\let\next@\msafam@\else
 \ifnum#2=\tw@\let\next@\msbfam@\fi\fi
 \mathchardef#1="#3\next@#4#5}
\def\mathhexbox@#1#2#3{\relax
 \ifmmode\mathpalette{}{\m@th\mathchar"#1#2#3}%
 \else\leavevmode\hbox{$\m@th\mathchar"#1#2#3$}\fi}
\def\hexnumber@#1{\ifcase#1 0\or 1\or 2\or 3\or 4\or 5\or 6\or 7\or 8\or
 9\or A\or B\or C\or D\or E\or F\fi}

\font\teneufm=eufm10 \font\seveneufm=eufm7 \font\fiveeufm=eufm5
\newfam\eufmfam
\textfont\eufmfam=\teneufm \scriptfont\eufmfam=\seveneufm
\scriptscriptfont\eufmfam=\fiveeufm

\catcode`\@=\csname pre amssym.def at\endcsname

\newcounter{remark}
\newcommand{\supp}{{\mathit supp}\,}

\renewcommand{\div}{\mbox{div}}

\def  \12  {{\frac{1}{2}}}

\def\build#1_#2^#3{\mathrel{\mathop{\kern 0pt#1}\limits_{#2}^{#3}}}

\numberwithin{equation}{section}

\begin{document}

\title[EMHD blowup]{Blowup for the forced electron MHD}


\author [Mimi Dai]{Mimi Dai}

\address{Department of Mathematics, Statistics and Computer Science, University of Illinois at Chicago, Chicago, IL 60607, USA}
\email{mdai@uic.edu} 



\thanks{The author is partially supported by the NSF grant DMS--2308208 and Simons Foundation.}

\begin{abstract}
The electron magnetohydrodynamics (MHD) contains a highly nonlinear Hall term with an interesting structure. Exploring the Hall nonlinear structure, we investigate possible phenomena of finite time blow up for the electron MHD with a (non-rough) forcing.
When the magnetic field has zero horizontal components, the vertical component equation has a mixing feature with the mixer being the current flow. By constructing a magnetic field profile whose current density is approximately a hyperbolic flow near the origin, we show blowup develops in finite time. In another setting when the magnetic field is a shear type, the Hall term vanishes, and finite time blowup can be obtained for the forced electron MHD as well.

\bigskip


\hspace{0.02cm}CLASSIFICATION CODE: 35Q35, 76B03, 76D09, 76E25, 76W05.
\end{abstract}

\maketitle

\section{Introduction}



We consider the electron magnetohydrodynamics (MHD) system with an external forcing
\begin{equation}\label{emhd}
\begin{split}
B_t+ \nabla\times ((\nabla\times B)\times B)=&\ F,\\
\nabla\cdot B=&\ 0.
\end{split}
\end{equation}
System \eqref{emhd} is an approximated model of the full magnetohydrodynamics with Hall effect while the slow ion motion and resistivity are negligible, see \cite{BDS, Bis1, KC}. 
In \eqref{emhd}, $B$ denotes the magnetic field satisfying the Gauss law (divergence free).
The Hall term $\nabla\times ((\nabla\times B)\times B)$ can be written as
\[\nabla\times ((\nabla\times B)\times B)=B\cdot \nabla J-J\cdot \nabla B\]
with the current density $J=\nabla\times B$, consisting of a transport term $-J\cdot \nabla B$ and a stretching term $B\cdot \nabla J$. The understanding of the two nonlinear features is the key to understand the Hall effect in MHD.

Although Hall effect has been a classical topic relevant to rapid magnetic reconnection in plasma physics, see for instance \cite{BDS, Bis1, BPS}, many mysteries surrounding magnetic reconnection phenomena remain to be answered. On the mathematical side, since the rigorous derivation of the Hall MHD system of PDEs in \cite{ADFL}, contributions have been made to study well-posedness and ill-posedness issues by many researchers. We briefly mention a few recent results, without the intention to provide a complete list in the literature.  
Well-posedness is obtained for electron MHD with a uniform magnetic background in \cite{Dai-emhd-2d, JO3}; various ill-posedness phenomena have been explored in \cite{Dai-emhd-24, Dai-unique24, Dai-unique18, JO1, JO2}. Reader can find more related work from the references in the above papers.

In this paper, we aim to explore possible mechanism of blowup from the Hall nonlinear structure. Hyperbolic flow has been studied for 2D Euler equation to construct solutions with double exponential growth in \cite{KS}; it has been also used in \cite{HL, HL2} to construct finite-time blowup solutions for the 3D axisymmetric Euler equation. In the recent papers \cite{CMZ, CMZZ} the authors constructed finite-time blowup solutions for the forced 3D Euler and hypodissipative Navier-Stokes equations by using profiles being hyperbolic near the origin. In our previous work \cite{Dai-ams} 1D models are proposed for the electron MHD to gain insights of the nonlinear effects. Constructing an initial odd profile being hyperbolic near the origin, the solution develops singularity in finite time. 

Inspired by the aforementioned findings, we continue to explore the possibility of finite-time blowup for the original electron MHD coming from hyperbolic profile. We discover that, when the magnetic field has zero horizontal components, the vertical component equation has a mixing feature with the mixer being the current flow. Moreover, choosing the vertical component of the magnetic field appropriately, its current density is approximately a hyperbolic flow near the origin. Such observation is the key to construct the following finite-time blowup scenario. 

\begin{Theorem}\label{thm}
For any small constant $\varepsilon>0$, there exists a solution $B(t)$ of the electron MHD \eqref{emhd} with an external forcing $F(t)$ on $[0,1]$ such that 
\begin{itemize}
\item [(i)] $F(t)\in C^{1-\varepsilon}(\mathbb R^3)\cap L^2(\mathbb R^3)$ for $\forall t\in[0,1]$; 
\item [(ii)] $B(t)\in C^{2-\varepsilon}(\mathbb R^3)\cap L^2(\mathbb R^3)$ for $\forall t\in[0,1]$; 
\item [(iii)] $\lim_{t\to 1}\|\nabla \times J(t)\|_{L^\infty}=\infty$;
\item [(iv)] $\int_0^1 \|\nabla \times J(t)\|_{L^\infty} \, dt=\infty$.
\end{itemize}
\end{Theorem}
In particular, if the magnetic field profile is chosen to be in the form $B(x,y,z, t)=(0,0, B^3(x,y, z, t))$, the current density is given by
\[J(x,y,z, t)=\nabla\times B(x,y,z, t)=\left( \partial_y B^3, -\partial_x B^3,0\right).\]
Consequently the electron MHD 
\[\partial_tB-J\cdot\nabla B=-B\cdot \nabla J\]
for this profile in component form reads 
\begin{equation}\notag
\begin{split}
\partial_tB^1&=-B^3\partial_z J^1,\\
\partial_tB^2&=-B^3\partial_z J^2,\\
\partial_tB^3-J\cdot\nabla B^3&=0.
\end{split}
\end{equation}
We notice that the vertical component $B^3$ satisfies a transport equation with the drift velocity being the horizontal current density, while the horizontal components $B^1$ and $B^2$ satisfies rather simple equations. Of course, a closer look tells us $J\cdot\nabla B^3=0$ actually. To avoid such triviality,
the component $B^3$ can be chosen such that the current density $J\approx (-x,y,0)$ near the origin and the vertical component equation is modified to
\[\partial_tB^3-\bar J\cdot\nabla B^3=0\]
where $\bar J\approx (-x,y,0)$ is a modified current of $J$ but $\bar J\cdot\nabla B^3\neq 0$.
The proof of Theorem \ref{thm} in Section \ref{sec-proof} relies on an iterative construction of such profiles generating hyperbolic current flow in Section \ref{sec-main}.

Further exploring the Hall term, we find that a shear type of profile $B(x,y,z,t)=(0,0, B^3(x,t))$ which generates a shear current density $J(x,y,z,t)=\nabla\times B(x,y,z,t)=(0, -\partial_xB^3(x,t), 0)$ can also be used to construct a solution to the forced electron MHD which develops blowup in finite time.  Namely,

\begin{Theorem}\label{thm2}
For any small constant $\varepsilon>0$, there exists a solution $B(t)$ of the electron MHD \eqref{emhd} with an external forcing $F(t)\in C^{2-\frac32\varepsilon}(\mathbb R^3)\cap L^2(\mathbb R^3)$ on $[0,1]$ such that (ii), (iii) and (iv) of Theorem \ref{thm} are satisfied.
\end{Theorem}


A few aspects regarding the results are discussed below.

\begin{Remark}
\begin{itemize}
\item [(i)] Both hyperbolic current flow and shear current profile are relevant in magnetic reconnection phenomena from the physics point of view. For the general $2\frac12D$ electron MHD, in which hyperbolic current flow is a special type, see the physics articles \cite{CSZ, WH}. On the other hand, current sheets have been extensively studied by physicists, for instance in \cite{BPS, GKLH, GH, WHL}. 
\item [(ii)] In \cite{CWeng} the authors studied the 3D axisymmetric electron (and Hall) MHD without external forcing and showed that the solution develops finite-time singularity in high order Sobolev space assuming existence in such space. In the current paper, we do not need to assume existence since the proof is constructive type.
\item [(iii)] Item (iii) in Theorem \ref{thm} indicates that the current density flow $J$ (serving as the drift velocity in the transport part of the electron MHD) fails to be $C^1$ at time $t=1$. This is consistent with the classical theory of transport equation in \cite{DL}.
\item [(iv)] For Euler equation, the classical Beale-Kato-Majda continuation criterion says that a solution $u(t)$ blows up at a time $T$ if and only if
\[\int_0^T \|\nabla \times u(t)\|_{L^\infty} \, dt=\infty.\]
In contrast, such criterion for the electron MHD is not known in general. Nevertheless, we obtained the Beale-Kato-Majda type criterion for the 3D resistive electron and Hall MHD in \cite{Dai-Oh} in term of $\nabla J$ (which involves more derivatives than $\nabla\times J$). In some special contexts of electron MHD, the Beale-Kato-Majda type criterion in term of $\nabla\times J$ was established in \cite{DGW, Dai-W}.
\item [(v)] We expect that blowup may even occur for the forced electron MHD with the presence of resistivity $(-\Delta)^\gamma B$ for some $\gamma>0$. This will be addressed in future project.
\item [(vi)] The ultimate goal is to construct finite-time blowup solutions for the electron MHD without external forcing. It will be further explored in forthcoming work.
\end{itemize}
\end{Remark}

\medskip


\textbf{Notations.}
As usual we denote $C$ by an insignificant general constant which may have different values at different lines.  The symbol $a\lesssim b$ denotes an inequality $a\leq Cb$ up to a constant when the constant $C$ does not play any important role. 

Sometimes we write $\vec x=(x,y,z)$. For convenience we also denote 
\[\partial_1=\partial_x, \ \  \partial_2=\partial_y, \ \ \partial_3=\partial_z.\]

\bigskip

\section{Main construction using hyperbolic flow}
\label{sec-main}
The goal of this section is to construct iteratively a sequence of approximating solutions to the electron MHD with external forcing. We shall start from a stationary magnetic field profile whose current density is close to a hyperbolic flow near the origin in the horizontal plane. Although more generic profile satisfying this condition would work, to reduce the complexity of estimates, the profile is chosen to have zero horizontal components and its non-zero vertical component does not depend on the $z$ variable. 
At each step, we add to the previous solution by a high frequency perturbation which has similar properties as the starting profile. 

\medskip

\subsection{Profile with a hyperbolic current density near the origin}
\label{sec-profile}
Let $f$ be a smooth cutoff function satisfying 
$f(x)=1$ for $|x|\leq \frac12$ and $f(x)=0$ for $|x|\geq 1$. We choose the stationary profile $B_0=(0,0, B_0^3)$ with the vertical component 
\begin{equation}\label{b0-3}
B_0^3=-\lambda^\gamma \sin(\lambda x)f(\lambda^{\gamma_1}x)\sin(\lambda y)f(\lambda^{\gamma_2}y)
\end{equation}
for a large constant $\lambda>0$ and parameters $\gamma\in \mathbb R, \gamma_1, \gamma_2>0$ to be determined later. It is clear that $\div B_0=0$ since $B_0^3$ does not depend on $z$ variable.

The corresponding current density is 
\[J_0=\nabla\times B_0=(\partial_y B_0^3, -\partial_x B_0^3, 0)\]
with 
\begin{equation}\label{j0}
\begin{split}
J_0^1=&\ \partial_y B_0^3=-\lambda^{1+\gamma}\sin(\lambda x)f(\lambda^{\gamma_1}x)\cos(\lambda y)f(\lambda^{\gamma_2}y)\\
&-\lambda^{\gamma+\gamma_2}\sin(\lambda x)f(\lambda^{\gamma_1}x)\sin(\lambda y)f'(\lambda^{\gamma_2}y)\\
J_0^2=&-\partial_x B_0^3= \lambda^{1+\gamma}\cos(\lambda x)f(\lambda^{\gamma_1}x)\sin(\lambda y)f(\lambda^{\gamma_2}y)\\
&+\lambda^{\gamma+\gamma_1}\sin(\lambda x)f'(\lambda^{\gamma_1}x)\sin(\lambda y)f(\lambda^{\gamma_2}y)\\
J_0^3=&\ 0.
\end{split}
\end{equation}
In view of \eqref{j0} we note for $|x|+|y|\ll \lambda^{-1}$, the current density satisfies 
\[J_0 \approx \lambda^{2+\gamma}(-x,y, 0)\]
which is a hyperbolic flow in the horizontal plane near the origin. This is the crucial point to generate blowup mechanism. 

It is easy to verify that 
\[J_0\cdot\nabla B_0=0, \ \ \ B_0\cdot\nabla J_0=0.\]
Hence $B_0$ is a stationary solution to the electron MHD \eqref{emhd} with zero external forcing $F_0=0$.

Define $\bar J_0=(\bar J_0^1, \bar J_0^2, \bar J_0^3)$ with
\begin{equation}\notag
\begin{cases}
\bar J_0^1=x\partial_x J_0^1(0,y,0),\\
\bar J_0^2= J_0^2(0,y,0),\\
\bar J_0^3=z\partial_z J_0^3(0,y,0).
\end{cases}
\end{equation}
One can see that 
\[\nabla\cdot \bar J_0=\partial_x\bar J_0^1+\partial_y\bar J_0^2+\partial_z\bar J_0^3=\partial_x J_0^1(0,y,0)+\partial_y J_0^2(0,y,0)+\partial_z J_0^3(0,y,0)=0\]
since $\nabla\cdot J_0=0$. On the other hand, straightforward computation shows that, using \eqref{j0}
\begin{equation}\label{bar-j0}
\begin{split}
\bar J_0^1=&-\lambda^{2+\gamma}x\cos(\lambda y)f(\lambda^{\gamma_2}y)
-\lambda^{1+\gamma+\gamma_2}x\sin(\lambda y)f'(\lambda^{\gamma_2}y),\\
\bar J_0^2=&\ \lambda^{1+\gamma}\sin(\lambda y)f(\lambda^{\gamma_2}y),\\
\bar J_0^3=&\ 0.
\end{split}
\end{equation}

\medskip

\subsection{Construction of a new solution $B_1$}
In this part, we show how to construct a new solution $B_1$ by adding to $B_0$ with a highly oscillatory perturbation $W_1$. In particular,  $W_1$ is produced from $\bar J_0$ which is a hyperbolic flow near the origin.

For a small constant $\varepsilon>0$, we define the frequency number
\begin{equation}\label{L-1}
\lambda_1=\exp\left\{2\int_{1-\lambda^{-\varepsilon/2}}^1-\partial_x J_0^1(0,s)\, ds\right\}=\exp\{2\lambda^{2+\gamma-\varepsilon/2}\}
\end{equation}
and 
\begin{equation}\label{k0}
K_0(t)=e^{\int_t^1-\partial_1\bar J_0^1(0,s)\, ds}=e^{(1-t)\lambda^{2+\gamma}}.
\end{equation}
For $2+\gamma-\varepsilon/2>0$, we have $\lambda_1\gg \lambda$.

Let $\bar W_1$ be the solution to the linear equation 
\begin{equation}\label{eq-bar-b1}
\partial_t\bar W_1-\bar J_0\cdot\nabla \bar W_1=-\bar W_1\cdot\nabla \bar J_0
\end{equation}
on the time interval $[0,1]$ with the ending time data
\begin{equation}\label{data}
\bar W_{1,\mbox{end}}(\vec x):=\bar W_1(\vec x,1)=(0,0,\bar W_1^3(\vec x,1))
\end{equation}
where 
\begin{equation}\label{data-3}
\bar W_1^3(\vec x,1)=-\lambda_1^{\gamma}\sin(\lambda_1 x)f(\lambda_1^{\gamma_1}x)\sin(\lambda_1 y)f(\lambda_1^{\gamma_2}y).
\end{equation}
We observe
\[\supp \bar W_1(\vec x,1)\subset \left\{\vec x=(x,y,z): |x|\leq \lambda_1^{-\gamma_1}, |y|\leq \lambda_1^{-\gamma_2}\right\}.\]

Since $\bar J_0^2$ is independent of $x$ and $z$, we have
\[\bar W_1\cdot\nabla \bar J_0^2=\bar W_1^1\partial_x \bar J_0^2+\bar W_1^2\partial_y \bar J_0^2+\bar W_1^3\partial_z \bar J_0^2=\bar W_1^2\partial_y \bar J_0^2.\]
It then follows from equation of $\bar W_1^2$ in \eqref{eq-bar-b1} and the zero ending data in \eqref{data} that $\bar W_1^2(x,t)\equiv 0$. Consequently we have
\[\bar W_1\cdot\nabla \bar J_0^1=\bar W_1^1\partial_x \bar J_0^1+\bar W_1^2\partial_y \bar J_0^1+\bar W_1^3\partial_z \bar J_0^1=\bar W_1^1\partial_x \bar J_0^1.\]
Moreover, $\bar W_1\cdot\nabla \bar J_0^3=0=\bar W_1^3\partial_z \bar J_0^3$, since $\bar J_0^3=0$. Therefore, \eqref{eq-bar-b1} is equivalent to 
\begin{equation}\label{bar-i}
\partial_t\bar W_1^i-\bar J_0\cdot\nabla \bar W_1^i=-\bar W_1^i\partial_i \bar J_0^i, \ \ i=1, 2,3
\end{equation}
where the right hand side is not a summation but a single term. Therefore, each component $\bar W_1^i$ satisfies a linear equation. 

Define the flow map $\phi_0(\vec x, t_0, t)=\left( \phi_0^1(\vec x, t_0, t),\phi_0^2(\vec x, t_0, t),\phi_0^3(\vec x, t_0, t)\right)$ as the solution to
\begin{equation}\label{flow}
\begin{cases}
\partial_t\phi_0^i(\vec x, t_0, t)=-\bar J_0^i(\phi_0(\vec x, t_0, t), t), \ \ i=1,2,3\\
\phi_0(\vec x, t_0, t_0)= \vec x.
\end{cases}
\end{equation} 

\begin{Lemma}\label{le-solution}
The solution of \eqref{eq-bar-b1} (and \eqref{bar-i}) with the ending time data \eqref{data}-\eqref{data-3} is given by
\begin{equation}\label{pert}
\bar W_1^i(\vec x, t)=\bar W_{1,\mbox{end}}^i\left( \phi_0(\vec x, t,1)\right) \exp\left\{ \int_1^t-\partial_i \bar J_0^i\left(\phi_0(\vec x, t, s),s \right) \, ds\right\}
\end{equation}
for $i=1,2,3$.
Moreover, we have
\begin{equation}\notag
\supp\left( \bar W_1(\vec x,t)\right)\subset\left\{\vec x=(x,y,z): |y|\leq 2\lambda_1^{-\gamma_2}K_0(t) \right\}
\end{equation}
with $K_0$ defined in \eqref{k0}.
\end{Lemma}
\pf
This formula can be found in \cite{CMZ}. We include a proof here for completeness. Denote 
\[\tilde W_1^i(x,t)=\bar W_1^i(\phi_0(\vec x, 1, t),t), \ \ i=1,2,3.\]
Then we have by using \eqref{eq-bar-b1} and \eqref{data}-\eqref{bar-b13}
\begin{equation}\notag
\begin{split}
\partial_t \tilde W_1^i(x,t)&=\partial_t \bar W_1^i(\phi_0(\vec x, 1, t),t)+\partial_t\phi_0(\vec x, 1, t) \cdot\nabla \bar W_1^i(\phi_0(\vec x, 1, t),t)\\
&=-\bar W_1^i(\phi_0(\vec x, 1, t),t)\partial_i\bar J_0^i(\phi_0(\vec x, 1, t),t)\\
&=-\tilde W_1^i(x,t)\partial_i\bar J_0^i(\phi_0(\vec x, 1, t),t)
\end{split}
\end{equation}
and 
\[\tilde W_1^i(x,1)=\bar W_1^i(\phi_0(\vec x, 1, 1),1)=\bar W_1^i(\vec x, 1)=\bar W_{1,\mbox{end}}^i(\vec x).\]
It follows that 
\begin{equation}\label{sol-1}
 \bar W_1^i(\phi_0(\vec x, 1, t),t)=\tilde W_1^i(x,t)=\bar W_1^i(\vec x, 1)\exp\left\{\int_1^t-\partial_i \bar J_0^i(\phi_0(\vec x,1,s), s) \,ds \right\}.
\end{equation}

Changing variable, we let $\phi_0(\vec x,1,t)=\vec y$. Acting $\phi_0$ on both sides yields
\[\phi_0(\phi_0(\vec x,1,t),t,1)=\phi_0(\vec y,t,1)\]
and hence $\vec x=\phi_0(\vec y,t,1)$. Applying the change of variables in \eqref{sol-1} we obtain
\begin{equation}\notag
\begin{split}
 \bar W_1^i(\vec y,t)&=\bar W_1^i(\phi_0(\vec y,t,1), 1)\exp\left\{\int_1^t-\partial_i \bar J_0^i(\phi_0(\phi_0(\vec y,t,1),1,s), s) \,ds \right\}\\
 &=\bar W_1^i(\phi_0(\vec y,t,1), 1)\exp\left\{\int_1^t-\partial_i \bar J_0^i(\phi_0(\vec y,t,s), s) \,ds \right\}.
 \end{split}
\end{equation}
The formula \eqref{pert} follows by replacing $\vec y$ with $\vec x$ in the formula above.

In view of \eqref{flow}, $\bar J_0^2(0)=0$ and $\bar J_0^2$ only depends on $y$, we have
\begin{equation}\notag
\begin{split}
\partial_t\phi_0^2(\vec x, 1,t)=&-\bar J_0^2(\phi_0(\vec x, 1,t), t)\\
=&-\phi_0^2(\vec x, 1,t)\partial_2 \bar J_0^2(0, t)-\bar J_{0,\mbox{rem}}^2(\phi_0(\vec x, 1,t), t)\\
=&\ \phi_0^2(\vec x, 1,t)\partial_1 \bar J_0^1(0, t)-\bar J_{0,\mbox{rem}}^2(\phi_0(\vec x, 1,t), t)\\
\end{split}
\end{equation}
with $\bar J_{0,\mbox{rem}}^2(\vec x, t)=\bar J_0^2(\vec x, t)-y\partial_2 \bar J_0^2(0, t)$, where we also used the fact $\div \bar J_0=0$. On the other hand, we see from \eqref{bar-j0} 
\[\partial_2^2 \bar J_0^2(0, t)=0, \ \ \left|\partial_2^3 \bar J_0^2(0, t)\right|\leq C\lambda^{4+\gamma}\]
and hence 
\begin{equation}\notag
|\bar J_{0,\mbox{rem}}^2(\vec x, t)|\leq C \left|\partial_2^3 \bar J_0^2(0, t)\right| |y|^3\leq C\lambda^{4+\gamma}|y|^3.
\end{equation}
It follows that for $t\in [t_0, 1]$ with $t_0\geq 1-\lambda^{-\frac12\varepsilon}$
\begin{equation}\label{phi-1}
\begin{split}
|\phi_0^2(\vec x, 1,t_0)|\leq&\ |y|e^{\int_1^{t_0}\partial_1 \bar J_0^1(0, t)-C\lambda^{4+\gamma}|\phi_0^2(\vec x, 1,t)|^2\, dt}\\
=&\ |y|e^{\int_{t_0}^1-\partial_1 \bar J_0^1(0, t)+C\lambda^{4+\gamma}|\phi_0^2(\vec x, 1,t)|^2\, dt}\\
=&\ |y|e^{(1-t_0)\lambda^{2+\gamma}}e^{\int_{t_0}^1C\lambda^{4+\gamma}|\phi_0^2(\vec x, 1,t)|^2\, dt}\\
\leq & |y| K_0(t_0)e^{C\lambda^{4+\gamma-\frac12\varepsilon}|\phi_0^2(\vec x, 1,t)|^2}.
\end{split}
\end{equation}
If $|\phi_0^2(\vec x, 1,t)|\leq \lambda^{-2-\frac12\gamma}$ for $t\in [t_0,1]$, \eqref{phi-1} implies for large enough $\lambda>0$
\begin{equation}\label{phi-2}
|\phi_0^2(\vec x, 1,t_0)|
\leq  |y| K_0(t_0)e^{C\lambda^{-\frac12\varepsilon}}\leq 2\lambda_1^{-\gamma_2} K_0(t_0)
\end{equation}
where we used $|y|\leq \lambda_1^{-\gamma_2}$ on which $\bar J_0^1$ is supported. On the other hand, recalling the definition of $\lambda_1$ in \eqref{L-1}, it follows from \eqref{phi-2} that for $t\in [t_0, 1]$
\begin{equation}\notag
|\phi_0^2(\vec x, 1,t)|\leq 2\lambda_1^{-\gamma_2} K_0(t)
\leq 2e^{-2\gamma_2\lambda^{2+\gamma-\frac12\varepsilon}} e^{(1-t)\lambda^{2+\gamma}}
\leq 2e^{(1-2\gamma_2)\lambda^{2+\gamma-\frac12\varepsilon}}.
\end{equation}
For $\gamma_2>\frac12$, we further claim that
\begin{equation}\notag
|\phi_0^2(\vec x, 1,t)|
\leq \lambda^{-(2+\gamma-\frac12\varepsilon)}\leq \lambda^{-2-\frac12\gamma}
\end{equation}
for large enough $\lambda>1$ and small enough $\varepsilon>0$, where we used the fact $c_1e^{-c_2\Lambda}\leq \Lambda^{-1}$ for $c_1,c_2>0$ and large enough $\Lambda>1$. Therefore we close the estimate
\begin{equation}\notag
|\phi_0^2(\vec x, 1,t)|
\leq 2 |y| K_0(t), \ \ t\in[t_0, 1], \ \ t_0\geq 1-\lambda^{-\frac12\varepsilon}.
\end{equation}
Consequently we have
\begin{equation}\notag
\supp\left( \bar W_1(\vec x,t)\right)\subset\left\{\vec x: |\phi_0^2(\vec x,t,1)|\leq \lambda_1^{-\gamma_2} \right\}\subset\left\{\vec x=(x,y,z): |y|\leq 2\lambda_1^{-\gamma_2}K_0(t) \right\}.
\end{equation}

\cbdu

As reflected in \eqref{data}, since $\bar W_{1,\mbox{end}}^1(\vec x)=\bar W_{1,\mbox{end}}^2(\vec x)=0$, it follows from Lemma \ref{le-solution} that
\[\bar W_{1}^1(\vec x,t)\equiv 0, \ \ \bar W_{1}^2(\vec x,t)\equiv 0.\]
On the other hand, since $\partial_3\bar J_0^3=0$, Lemma \ref{le-solution} implies
\[\bar W_{1}^3(\vec x,t)=\bar W_{1,\mbox{end}}^3(\phi_0(\vec x,t,1)).\]
Note $\phi_0^1(\vec x, t,1)$ and $\phi_0^2(\vec x, t,1)$ depend only on $x$ and $y$ but not on $z$. Hence $\bar W_1^3$ depends only on $x$ and $y$ but not on $z$. It is easy to see
\begin{equation}\notag
\begin{split}
\nabla\times \bar W_1&=(\partial_y \bar W_1^3, -\partial_x \bar W_1^3,0),\\
\nabla\times \bar W_1\cdot \bar W_1&=0,\\
\bar W_1\cdot\nabla(\nabla\times \bar W_1)&=0.
\end{split}
\end{equation}

Denote 
\begin{equation}\notag
\phi_0^1(\vec x,t,1)=xK_0(t) \tilde \phi_0^1(y,t).
\end{equation}
More analysis on the flow map shows
\begin{Lemma}\label{le-phi}
Let $\gamma_2>\frac12$. For $t\in[1-\lambda^{-\frac12\varepsilon}, 1]$ and $|y|\leq 2\lambda_1^{-\gamma_2}K_0(t)$, we have
\begin{equation}\notag
\begin{split}
\frac12\leq |\tilde \phi_0^1(y,t)|\leq 2, \ \ |\partial_2\tilde \phi_0^1(y,t)|\leq 8 \lambda^{4+\gamma}\lambda_1^{-\gamma_2}K_0(t), \ \
|\partial_2^2\tilde \phi_0^1(y,t)|\leq 8\lambda^{4+\gamma},\\
\frac{1}{2K_0(t)}\leq |\partial_2\phi_0^2(\vec x,t,1)|\leq\frac{2}{K_0(t)},\\
|\partial_2^2\phi_0^2(\vec x,t,1)|\leq \lambda^{4+\gamma}\lambda_1^{-\gamma_2}K_0(t).
\end{split}
\end{equation}
\end{Lemma}
\pf
In view of \eqref{flow}, using Taylor's expansion we deduce
\begin{equation}\notag
\begin{split}
&\exp\left\{ \int_{t_0}^t\partial_2(-\bar J_0^2)(0,s)-\frac12y^2\partial_2^3(-\bar J_0^2)(0,s)+\cdot\cdot\cdot \, ds\right\}\\
\leq& |\partial_0^2\phi_0^2(\vec x,t_0,t)|\\
\leq& \exp\left\{ \int_{t_0}^t\partial_2(-\bar J_0^2)(0,s)+\frac12y^2\partial_2^3(-\bar J_0^2)(0,s)+\cdot\cdot\cdot \, ds\right\}.
\end{split}
\end{equation}
It follows
\begin{equation}\notag
\begin{split}
&\exp\left\{ \int_{t_0}^t\partial_2(-\bar J_0^2)(0,s)-y^2\lambda^{4+\gamma}\, ds\right\}\\
\leq& |\partial_2\phi_0^2(\vec x,t_0,t)|\\
\leq& \exp\left\{ \int_{t_0}^t\partial_2(-\bar J_0^2)(0,s)+y^2\lambda^{4+\gamma} \, ds\right\},
\end{split}
\end{equation}
and thus
\begin{equation}\notag
\begin{split}
&\exp\left\{ \int_{t_0}^t-\lambda^{2+\gamma}-y^2\lambda^{4+\gamma}\, ds\right\}\\
\leq& |\partial_2\phi_0^2(\vec x,t_0,t)|\\
\leq& \exp\left\{ \int_{t_0}^t-\lambda^{2+\gamma}+y^2\lambda^{4+\gamma} \, ds\right\}.
\end{split}
\end{equation}
For $|y|\leq 2\lambda_1^{-\gamma_2}K_0(t)$, it has
\[y^2\lambda^{4+\gamma}\leq 4\lambda_1^{-2\gamma_2}K_0^2(t)\lambda^{4+\gamma}\]
and thus
\[-\lambda^{2+\gamma}+y^2\lambda^{4+\gamma}\leq -\lambda^{2+\gamma}(1-4\lambda_1^{-2\gamma_2}K_0^2(t)\lambda^2).\]
Recall $\lambda_1=e^{2\lambda^{+2+\gamma-\frac12\varepsilon}}$ and 
\[K_0(t)=e^{(1-t)\lambda^{2+\gamma}}=e^{2\lambda^{2+\gamma-\frac12\varepsilon}\cdot \frac12(1-t)\lambda^{\frac12\varepsilon}}=\lambda_1^{\frac12(1-t)\lambda^{\frac12\varepsilon}}\]
which implies
\[\lambda_1^{-2\gamma_2}K_0^2(t)\lambda^2=e^{-4\gamma_2\lambda^{-\varepsilon/2+2+\gamma}+2(1-t)\lambda^{2+\gamma}}\lambda^2<\frac14\]
for $\gamma_2>\frac12$.

Hence for $\gamma_2>\frac12$, we have
\[|\partial_2\phi_0^2(\vec x,t_0,t)|\leq 1\]
\[|\partial_2^2\phi_0^2(\vec x,t_0,t)|\leq |t-t_0|\lambda^{4+\gamma}|y|,\]
\[\frac{1}{2K_0(t)}\leq |\partial_2\phi_0^2(\vec x,t,1)|\leq\frac{2}{K_0(t)},\]
\[|\partial_2^2\phi_0^2(\vec x,t,1)|\leq \lambda^{4+\gamma}\lambda_1^{-\gamma_2}K_0(t)=\lambda^{4+\gamma}e^{-2\gamma_2\lambda^{-\varepsilon/2+2+\gamma}+(1-t)\lambda^{2+\gamma}}.\]

Regarding the component $\phi_0^1(\vec x,t,1)$, the equation
\begin{equation}\notag
\partial_t \phi_0^1(\vec x,t,1)=-\bar J_0^1 (\phi_0(\vec x,t,1),t)=-\phi_0^1(\vec x,t,1)(\partial_1\bar J_0^1)\left((0, \phi_0^2(\vec x,t,1),0),t \right),
\end{equation}
gives
\begin{equation}\notag
\phi_0^1(\vec x,t,1)=xK_0(t) \exp\left\{ \int_t^1 P_{0,1}(\phi_0(y,t,s),s) \, ds\right\}
\end{equation}
with 
\begin{equation}\notag
P_{0,1}(y,t)=\partial_1(-\bar J_0^1)\left((0, \phi_0^2(y,t,1),0),t \right)-\partial_1(-\bar J_0^1)(0,t).
\end{equation}
Denote 
\begin{equation}\notag
\tilde \phi_0^{1}(y,t):=\exp\left\{ \int_t^1 P_{0,1}(\phi_0(y,t,s),s) \, ds\right\}.
\end{equation}
We have the estimates
\begin{equation}\notag
\begin{split}
|P_{0,1}(\phi_0(y,t_0,t),t)|\leq \lambda^{4+\gamma}|\phi_0(y,t_0,t)|^2\lesssim \lambda^{4+\gamma}y^2,\\
|\tilde \phi_0^{1}(y,t)|\leq e^{\int_t^1\lambda^{4+\gamma}y^2\, ds}\in [\frac12,2],\\
|\partial_2\tilde \phi_0^{1}(y,t)|\leq 8\lambda^{4+\gamma}\lambda_1^{-\gamma_2}K_0(t),\\
|\partial_2^2\tilde \phi_0^{1}(y,t)|\leq 8\lambda^{4+\gamma}.
\end{split}
\end{equation}


\cbdu

For the perturbation $\bar W_1$ we obtain the following estimate for the H\"older norm.

\begin{Lemma}\label{le-b1-holder}
Let $k\in \mathbb N$ and $\beta\in(0,1)$. It holds
\[ \|\bar W_1(t)\|_{C^{k,\beta}}\lesssim \lambda_1^{\gamma+k+\beta}, \ \ \mbox{for} \ 1-\lambda_1^{-\frac12\varepsilon}\leq t\leq 1.\]
\end{Lemma}
\pf
Since $\bar W_1^1(t)\equiv 0$ and $\bar W_2^1(t)\equiv 0$, we only need to show the estimate for $\bar W_1^3(t)$. It follows from
\[\bar W_{1}^3(\vec x,t)=\bar W_{1,\mbox{end}}^3(\phi_0(\vec x,t,1))\]
that
\begin{equation}\notag
\|\bar W_1^3(t)\|_{C^{k,\beta}}\lesssim \|\bar W_{1,\mbox{end}}^3(\phi_0(\vec x,1,1))\|_{C^{k,\beta}} \left(1+ \|\partial_x \phi_0(\vec x,1,t) \|_{C^k} + \|\partial_y \phi_0(\vec x,1,t) \|_{C^k}\right)^{k+1}.
\end{equation}
In view of \eqref{data-3}, we have
\begin{equation}\notag
\|\bar W_{1,\mbox{end}}^3(\phi_0(\vec x,1,1))\|_{C^{k,\beta}}=\|\bar W_{1,\mbox{end}}^3(\vec x)\|_{C^{k,\beta}} \lesssim \lambda_1^{\gamma+k+\beta}.
\end{equation}
It follows from \eqref{flow} that
\begin{equation}\notag
\partial_t\partial_x\phi_0(\vec x, 1, t)=-\partial_x\bar J_0 \cdot \partial_x(\phi_0(\vec x, 1, t), t)
\end{equation}
and hence
\begin{equation}\notag
\begin{split}
\partial_t\|\partial_x \phi_0(\vec x,1,t)\|_{L^\infty}\lesssim \|\partial_x \bar J_0\|_{L^\infty} \|\partial_x \phi_0(\vec x,1,t)\|_{L^\infty}
\lesssim \lambda^{\gamma+2} \|\partial_x \phi_0(\vec x,1,t)\|_{L^\infty}
\end{split}
\end{equation}
Integrating in time yields 
\[\|\partial_x \phi_0(\vec x,1,t)\|_{L^\infty}\lesssim 1+e^{\lambda_1^{-\frac12\varepsilon}\lambda^{\gamma+2}}\lesssim 1, \ \ t\in [1-\lambda_1^{-\frac12\varepsilon}, 1].\]
Similarly we can show such an estimate for $\partial_y \phi_0(\vec x,1,t)$. Iteratively we can obtain 
\[\|\partial_x \phi_0(\vec x,1,t) \|_{C^k}+\|\partial_y \phi_0(\vec x,1,t) \|_{C^k}\lesssim 1, \ \ \mbox{for any} \ \ k\in \mathbb N. \]
Therefore, for $t\in [1-\lambda_1^{-\frac12\varepsilon}, 1]$
\begin{equation}\notag
\|\bar W_1^3(t)\|_{C^{k,\beta}}\lesssim \|\bar W_{1,\mbox{end}}^3(\phi_0(\vec x,1,1))\|_{C^{k,\beta}} \lesssim \lambda_1^{\gamma+k+\beta}.
\end{equation}

\cbdu

We are ready to define the final perturbation 
\[W_1(\vec x, t)=g_1(t)\bar W_1(\vec x,t)\]
with the smooth temporal cutoff function $g_1$ defined by
\begin{equation}\label{g1}
g_1(t)=
\begin{cases}
0, \ \ \ \ t\leq 1-\lambda^{-\varepsilon/2}\\
(t-1+\lambda^{-\varepsilon/2})\lambda_{1}^{\varepsilon}, \ \ \ 1-\lambda^{-\varepsilon/2}\leq t\leq 1-\lambda^{-\varepsilon/2}+\lambda_1^{-\varepsilon}\\
1, \ \ \ \ 1-\lambda^{-\varepsilon/2}+\lambda_1^{-\varepsilon}\leq t\leq 1.
\end{cases}
\end{equation}
Note $W_1(\vec x,t)=0$ for $t\leq 1-\lambda^{-\varepsilon/2}$.
Denote $V_1=\nabla \times W_1=g_1(t)\nabla\times \bar W_1(\vec x,t)$. 


We define the new solution 
\begin{equation}\label{sol-B1}
B_1(\vec x, t)=B_0(\vec x)+W_1(\vec x, t)=(0,0, B_0^3(\vec x)+g_1(t) \bar W_1^3(\vec x, t))
\end{equation}
with 
\[\bar W_{1}^3(\vec x,t)=\bar W_{1,\mbox{end}}^3(\phi_0(\vec x,t,1)).\]
We point out that $\bar W_{1}^3$ does not depend on the vertical variable $z$. Denote the current $J_1=\nabla\times B_1$. 
The new solution $B_1$ satisfies the electron MHD \eqref{emhd} with an external forcing $F_1$,
\begin{equation}\label{iter-b1}
\partial_tB_1-J_1\cdot \nabla B_1=-B_1\cdot\nabla J_1+F_1
\end{equation}
where 
\begin{equation}\notag
\begin{split}
F_1=&\ g_1'(t)\bar W_1+g_1(t)\left[(\bar J_0-J_0)\cdot\nabla \bar W_1-\bar W_1\cdot\nabla (\bar J_0-J_0) \right]\\
&+g_1(t)\left[-(\nabla\times \bar W_1)\cdot\nabla B_0+B_0\cdot\nabla(\nabla\times \bar W_1) \right]\\
&+g_1^2(t)\left[-(\nabla\times \bar W_1)\cdot\nabla \bar W_1+\bar W_1\cdot\nabla(\nabla\times \bar W_1) \right].
\end{split}
\end{equation}
Recall 
\[(\nabla\times \bar W_1)\cdot\nabla \bar W_1=0, \ \ \bar W_1\cdot\nabla(\nabla\times \bar W_1)=0,\]
and also 
\[ B_0\cdot\nabla(\nabla\times \bar W_1)=0, \ \ \bar W_1\cdot\nabla(\bar J_0-J_0)=0.\]
Hence the force $F_1$ reduces to
\begin{equation}\label{F1}
F_1=g_1'(t)\bar W_1+g_1(t)(\bar J_0-J_0)\cdot\nabla \bar W_1-g_1(t)(\nabla\times \bar W_1)\cdot\nabla B_0.
\end{equation}

\begin{Lemma}\label{le-F1-holder}
For $1-\lambda_1^{-\frac12\varepsilon}\leq t\leq 1$, it holds
\[ \|F_1(t)\|_{C^\alpha}\lesssim \lambda_1^{1+\gamma+\alpha}.\]
\end{Lemma}
\pf
Recall
\[g_1'(t)\sim \lambda_1^{\varepsilon}, \ \ \|\bar W_1\|_{C^\alpha}\lesssim \lambda_1^{\gamma+\alpha}, \ \ \|\bar B_0\|_{C^\alpha}\lesssim \lambda^{\gamma+\alpha}.\]
It thus implies 
\[|g_1'(t)\bar W_1|\lesssim \lambda_1^{\varepsilon+\gamma}, \ \ |g_1(t)(\nabla\times \bar W_1)\cdot\nabla B_0| \lesssim \lambda^{\gamma+1}\lambda_1^{\gamma+1}.\]
Additionally, using \eqref{j0} and \eqref{bar-j0} we infer
\begin{equation}\notag
\begin{split}
\bar J_0^1-J_0^1=&-\lambda^{2+\gamma}x\cos(\lambda y)f(\lambda^{\gamma_2}y)
-\lambda^{1+\gamma+\gamma_2}x\sin(\lambda y)f'(\lambda^{\gamma_2}y)\\
&+\lambda^{1+\gamma}\sin(\lambda x)f(\lambda^{\gamma_1}x)\cos(\lambda y)f(\lambda^{\gamma_2}y)\\
&+\lambda^{\gamma+\gamma_2}\sin(\lambda x)f(\lambda^{\gamma_1}x)\sin(\lambda y)f'(\lambda^{\gamma_2}y)\\
=&\ \lambda^{1+\gamma}\left[\sin(\lambda x)f(\lambda^{\gamma_1}x)-\lambda x \right]\cos(\lambda y)f(\lambda^{\gamma_2}y)\\
&+\lambda^{\gamma+\gamma_2}\left[\sin(\lambda x)f(\lambda^{\gamma_1}x)-\lambda x \right]\sin(\lambda y)f'(\lambda^{\gamma_2}y),
\end{split}
\end{equation}
\begin{equation}\notag
\begin{split}
\bar J_0^2-J_0^2=&\ \lambda^{1+\gamma}\sin(\lambda y)f(\lambda^{\gamma_2}y)\\
&-\lambda^{1+\gamma}\cos(\lambda x)f(\lambda^{\gamma_1}x)\sin(\lambda y)f(\lambda^{\gamma_2}y)\\
&-\lambda^{\gamma+\gamma_1}\sin(\lambda x)f'(\lambda^{\gamma_1}x)\sin(\lambda y)f(\lambda^{\gamma_2}y)\\
=&\ \lambda^{1+\gamma}\sin(\lambda y)f(\lambda^{\gamma_2}y)\left[1-\cos(\lambda x)f(\lambda^{\gamma_1}x) \right]\\
&-\lambda^{\gamma+\gamma_1}\sin(\lambda x)f'(\lambda^{\gamma_1}x)\sin(\lambda y)f(\lambda^{\gamma_2}y).
\end{split}
\end{equation}
As a consequence, we note
\begin{equation}\notag
\begin{split}
\bar J_0^1-J_0^1\approx&\ \lambda^{1+\gamma}\lambda^3 x^3\cos(\lambda y)f(\lambda^{\gamma_2}y)
+\lambda^{\gamma+\gamma_2}\lambda^3 x^3\sin(\lambda y)f'(\lambda^{\gamma_2}y),\\
\bar J_0^2-J_0^2\approx&\ \lambda^{1+\gamma}\cdot\lambda y\cdot \lambda^2 x^2-\lambda^{\gamma+\gamma_1}\cdot\lambda x\cdot \lambda y.
\end{split}
\end{equation}
Hence we have 
\begin{equation}\notag
\begin{split}
|\bar J_0^1-J_0^1|\lesssim&\ \lambda^{1+\gamma}\lambda^3|x|^3+ \lambda^{\gamma+\gamma_2}\lambda^3|x|^3\lambda |y|
\lesssim \lambda^{1+\gamma}\lambda^3|x|^3,\\
|\bar J_0^2-J_0^2|\lesssim &\ \lambda^{1+\gamma}\lambda |y|\lambda^2 |x|^2+\lambda^{\gamma+\gamma_1}\lambda |x|\lambda |y|
\lesssim \lambda^{1+\gamma}\lambda^2 |x| |y|.
\end{split}
\end{equation}

For $|x|\leq \lambda_1^{-\gamma_1}$ and $|y|\leq 2\lambda_1^{-\gamma_2}K_0(t)$, we have
\begin{equation}\notag
\begin{split}
|g_1(t)(\bar J_0-J_0)\cdot\nabla \bar W_1|&\lesssim |\bar J_0-J_0||\nabla \bar W_1|\\
&\lesssim \lambda^{4+\gamma}|x|(|x|^2+|y|)\lambda_1^{1+\gamma}\\
&\lesssim \lambda^{4+\gamma}\lambda_1^{-\gamma_1}\left( \lambda_1^{-2\gamma_1}+\lambda_1^{-\gamma_2}e^{(1-t)\lambda^{2+\gamma}}\right)\lambda_1^{1+\gamma}\\
&\lesssim \lambda^{4+\gamma}\lambda_1^{-\gamma_1+1+\gamma}\left(\lambda_1^{-2\gamma_1}+\lambda_1^{-\gamma_2+\varepsilon} \right)\\
&\lesssim \lambda^{4+\gamma}\lambda_1^{-\gamma_1+1+\gamma-\gamma_2+\varepsilon}
\end{split}
\end{equation}
provided $2\gamma_1\geq \gamma_2-\varepsilon$. Combining the analysis above we infer
\begin{equation}\notag
\|F_1(t)\|_{C^{\alpha}}\lesssim \lambda_1^{\varepsilon+\gamma+\alpha}+\lambda_1^{1+\gamma+\alpha}+\lambda_1^{-\gamma_1+1+\gamma-\gamma_2+\varepsilon+\alpha}.
\end{equation}
Hence for $\alpha<-1-\gamma$, it follows
\[\|F_1(t)\|_{C^{\alpha}}\lesssim 1.\]
Recall $\gamma>-2$, we expect $F_1$ to be in $C^{1-\varepsilon}$.

\cbdu

\medskip

\subsection{Construction of the next solution $B_2$}
\label{sec-b2}
To illustrate the iteration in a transparent manner, we show the detail 
of constructing the next level solution $B_2$ by constructing the perturbation $W_2=B_2-B_1$. Ideally we would hope to produce $W_2$ from the previous step by a hyperbolic flow again in order to carry on the iteration. 
In particular, the new perturbation will be obtained through the previous perturbation 
\[W_1(\vec x,t)=(0,0, g_1(t)\bar W_{1,\mbox{end}}^3(\phi_0(\vec x, t,1)))\]
with
\begin{equation}\notag
\begin{split}
&\bar W_{1,\mbox{end}}^3\left(\phi_0(\vec x,t,1)\right)\\
=&-\lambda_1^{\gamma}\sin\left(\lambda_1\phi_0^1(\vec x,t,1) \right)f\left(\lambda_1^{\gamma_1}\phi_0^1(\vec x,t,1) \right)
\sin\left(\lambda_1\phi_0^2(\vec x,t,1) \right)f\left(\lambda_1^{\gamma_2}\phi_0^2(\vec x,t,1) \right)\\
=&-\lambda_1^{\gamma}\sin\left(\lambda_1x\tilde \phi_0^{1}(y)K_0(t) \right)f\left(\lambda_1^{\gamma_1}x\tilde \phi_0^{1}(y)K_0(t) \right)\sin(\lambda_1y\tilde\phi_0^2(y))f(\lambda_1^{\gamma_2}y\tilde\phi_0^2(y))\\
=&-\lambda_1^{\gamma}\sin\left(\lambda_1^{1+\frac12(1-t)\lambda^{\varepsilon/2}}x\tilde \phi_0^{1}(y)\right)f\left(\lambda_1^{\gamma_1+\frac12(1-t)\lambda^{\varepsilon/2}}x\tilde \phi_0^{1}(y) \right)\\
&\cdot\sin(\lambda_1y\tilde \phi_0^2(y))f(\lambda_1^{\gamma_2}y\tilde \phi_0^2(y))
\end{split}
\end{equation}
where
\[\tilde \phi_0^2(y)=e^{\int_1^t-J_0^2((0,y,0),s)\, ds}, \ \ \phi_0^2(\vec x,t,1)=y\tilde \phi_0^2(y).\]
As shown in Lemma \ref{le-phi}, $\frac12\leq |\tilde\phi_0^1(y)|\leq 2$, and hence $W_1^3(t)$ is compactly supported on 
\begin{equation}\label{supp}
\{\vec x=(x,y,z): |x|\leq 2\lambda_1^{-\gamma_1-\frac12(1-t)\lambda^{\varepsilon/2}}, |y|\leq 2\lambda_1^{-\gamma_2}K_0(t) \}. 
\end{equation}

We need to argue that the current density of the perturbation $W_1$ given by
\[V_1=\nabla\times W_1=(\partial_y W_1^3, -\partial_x W_1^3, 0)\]
is close to a hyperbolic flow near the origin. To do so, we compute the components of $V_1$ explicitly in the following,
\begin{equation}\label{V11}
\begin{split}
V_1^1=\partial_y W_1^3=&\ g_1(t)\partial_y \bar W_{1,\mbox{end}}^3\left(\phi_0(\vec x,t,1)\right)\\
=&-g_1(t)\lambda_1^{\gamma}\cos\left(\lambda_1^{1+\frac12(1-t)\lambda^{\varepsilon/2}}x\tilde \phi_0^{1}(y) \right)\lambda_1^{1+\frac12(1-t)\lambda^{\varepsilon/2}}x\partial_y\tilde \phi_0^{1}(y)\\
&\cdot f\left(\lambda_1^{\gamma_1+\frac12(1-t)\lambda^{\varepsilon/2}}x\tilde \phi_0^1(y) \right)
\sin(\lambda_1y\tilde \phi_0^2(y))f(\lambda_1^{\gamma_2}y\tilde \phi_0^2(y))\\
&-g_1(t)\lambda_1^{\gamma}\sin\left(\lambda_1^{1+\frac12(1-t)\lambda^{\varepsilon/2}}x\tilde \phi_0^{1}(y) \right)\lambda_1^{1+\frac12(1-t)\lambda^{\varepsilon/2}}x\partial_y\tilde \phi_0^{1}(y)\\
&\cdot f'\left(\lambda_1^{\gamma_1+\frac12(1-t)\lambda^{\varepsilon/2}}x\tilde \phi_0^{1}(y) \right)
\sin(\lambda_1y\tilde \phi_0^2(y))f(\lambda_1^{\gamma_2}y\tilde \phi_0^2(y))\\
&-g_1(t)\lambda_1^{\gamma}\sin\left(\lambda_1^{1+\frac12(1-t)\lambda^{\varepsilon/2}}x\tilde \phi_0^{1}(y)\right)f\left(\lambda_1^{\gamma_1+\frac12(1-t)\lambda^{\varepsilon/2}}x\tilde \phi_0^{1}(y) \right)\\
&\cdot\cos(\lambda_1y\tilde \phi_0^2(y))(\lambda_1\tilde \phi_0^2(y)+\lambda_1y\partial_y\tilde \phi_0^2(y))f(\lambda_1^{\gamma_2}y\tilde \phi_0^2(y))\\
&-g_1(t)\lambda_1^{\gamma}\sin\left(\lambda_1^{1+\frac12(1-t)\lambda^{\varepsilon/2}}x\tilde \phi_0^{1}(y)\right)f\left(\lambda_1^{\gamma_1+\frac12(1-t)\lambda^{\varepsilon/2}}x\tilde \phi_0^{1}(y) \right)\\
&\cdot\sin(\lambda_1y\tilde \phi_0^2(y))f'(\lambda_1^{\gamma_2}y\tilde \phi_0^2(y)) (\lambda_1^{\gamma_2}\tilde \phi_0^2(y)+\lambda_1^{\gamma_2}y\partial_y\tilde \phi_0^2(y)),
\end{split}
\end{equation}
and
\begin{equation}\label{V12}
\begin{split}
V_1^2=-\partial_xW_1^3=&-g_1(t)\partial_x\bar W^3_{1,\mbox{end}}\left(\phi_0(\vec x,t,1)\right)\\
=&\ g_1(t)\lambda_1^{\gamma+1+\frac12(1-t)\lambda^{\varepsilon/2}}\tilde \phi_0^1\cos(\lambda_1^{1+\frac12(1-t)\lambda^{\varepsilon/2}}x\tilde \phi_0^1)f(\lambda_1^{\gamma_1+\frac12(1-t)\lambda^{\varepsilon/2}}x\tilde \phi_0^1)\\
&\cdot\sin(\lambda_1\tilde \phi_0^2y)f(\lambda_1^{\gamma_2}\tilde \phi_0^2y)\\
&+g_1(t)\lambda_1^{\gamma+\gamma_1+\frac12(1-t)\lambda^{\varepsilon/2}}\tilde \phi_0^1\sin(\lambda_1^{1+\frac12(1-t)\lambda^{\varepsilon/2}}x\tilde \phi_0^1)f'(\lambda_1^{\gamma_1+\frac12(1-t)\lambda^{\varepsilon/2}}x\tilde \phi_0^1)\\
&\cdot\sin(\lambda_1\tilde \phi_0^2y)f(\lambda_1^{\gamma_2}\tilde \phi_0^2y).
\end{split}
\end{equation}

Define the modified current density
$\bar V_1=(\bar V_1^1, \bar V_1^2, \bar V_1^3)$ by
\begin{equation}\label{bar-V1}
\begin{cases}
\bar V_1^1=x\partial_x V_1^1((0,y,0),t),\\
\bar V_1^2= V_1^2((0,y,0),t),\\
\bar V_1^3=z\partial_z V_1^3((0,y,0),t).
\end{cases}
\end{equation}
Continuing with \eqref{V11}, we further compute
\begin{equation}\label{dx-V11}
\begin{split}
&\partial_x V_1^1((0,y,0),t)\\
=&-g_1(t)\lambda_1^{\gamma+1+\frac12(1-t)\lambda^{\varepsilon/2}}\partial_y \tilde \phi_0^1\sin(\lambda_1 y\tilde \phi_0^2)f(\lambda_1^{\gamma_2}y\tilde \phi_0^2)\\
&-g_1(t)\lambda_1^{\gamma+2+\frac12(1-t)\lambda^{\varepsilon/2}}\tilde \phi_0^1\left[y\partial_y \tilde \phi_0^2 +\tilde \phi_0^2\right]\cos(\lambda_1y \tilde \phi_0^2)f(\lambda_1^{\gamma_2}y\tilde \phi_0^2)\\
&-g_1(t)\lambda_1^{\gamma+\gamma_2+1+\frac12(1-t)\lambda^{\varepsilon/2}}\tilde \phi_0^1\left[y\partial_y \tilde \phi_0^2 +\tilde \phi_0^2\right]\sin(\lambda_1y \tilde \phi_0^2)f'(\lambda_1^{\gamma_2}y\tilde \phi_0^2).
\end{split}
\end{equation}
Hence we have
\begin{equation}\label{bar-V11}
\begin{split}
\bar V_1^1=&\ x\partial_x V_1^1((0,y,0),t)\\
=&-g_1(t)\lambda_1^{\gamma+1+\frac12(1-t)\lambda^{\varepsilon/2}}x\partial_y \tilde \phi_0^1\sin(\lambda_1 y\tilde \phi_0^2)f(\lambda_1^{\gamma_2}y\tilde \phi_0^2)\\
&-g_1(t)\lambda_1^{\gamma+2+\frac12(1-t)\lambda^{\varepsilon/2}}x\tilde \phi_0^1\left[y\partial_y \tilde \phi_0^2 +\tilde \phi_0^2\right]\cos(\lambda_1y \tilde \phi_0^2)f(\lambda_1^{\gamma_2}y\tilde \phi_0^2)\\
&-g_1(t)\lambda_1^{\gamma+\gamma_2+1+\frac12(1-t)\lambda^{\varepsilon/2}}x\tilde \phi_0^1\left[y\partial_y \tilde \phi_0^2 +\tilde \phi_0^2\right]\sin(\lambda_1y \tilde \phi_0^2)f'(\lambda_1^{\gamma_2}y\tilde \phi_0^2).
\end{split}
\end{equation}

It follows from \eqref{V12} that
\begin{equation}\label{bar-V12}
\bar V_1^2(\vec x,t)=V_1^2((0,y,0), t)=g_1(t)\lambda_1^{\gamma+1+\frac12(1-t)\lambda^{\varepsilon/2}}\sin(\lambda_1y\tilde \phi_0^2)f(\lambda_1^{\gamma_2}y\tilde \phi_0^2).
\end{equation}

By Lemma \ref{le-phi}, we know $\frac12\leq |\tilde \phi_0^1|\leq 2$ and
\[ |\partial_y \tilde \phi_0^1|\leq 8\lambda^{4+\gamma}\lambda_1^{-\gamma_2}K_0(t)=8\lambda^{4+\gamma}\lambda_1^{-\gamma_2+\frac12(1-t)\lambda^{\varepsilon/2}}\lesssim 1  \]
for $\gamma_2>\frac12$ and $1-\lambda^{-\frac12\varepsilon}\leq t\leq 1$. Moreover, 
\begin{equation}\label{phi-2}
|y\partial_y \tilde \phi_0^2 +\tilde \phi_0^2|=|\partial_y \phi_0^2|\in \left[\frac1{2K_0(t)}, \frac2{K_0(t)}\right].
\end{equation}
Therefore, for small enough $|x|$ and $|y|$ on the support \eqref{supp}, the second term in \eqref{bar-V11} is the dominating one and
\[\bar V_1^1\approx -g_1(t)\lambda_1^{\gamma+2+\frac12(1-t)\lambda^{\varepsilon/2}}x K_0(t)\approx -g_1(t)\lambda_1^{\gamma+2}x .\]
Again, for small enough $|y|$, it follows from \eqref{phi-2} that 
\begin{equation}\label{phi2-size}
|\tilde \phi_0^2|\in \left[\frac1{4K_0(t)}, \frac4{K_0(t)}\right]. 
\end{equation}
Hence, we also have
\[\bar V_1^2(\vec x,t)\approx g_1(t)\lambda_1^{\gamma+1+\frac12(1-t)\lambda^{\varepsilon/2}}\lambda_1y\tilde \phi_0^2\approx g_1(t)\lambda_1^{\gamma+2} y\]
in view of \eqref{bar-V12}. Therefore, we have again a hyperbolic flow in the horizontal plane near the origin
\[\bar V_1\approx g_1(t)\lambda_1^{\gamma+2} (-x,y,0) \]
which is the motivation of constructing the next perturbation using the modified current density $\bar V_1$.

Define the new frequency
\begin{equation}\notag
\lambda_2=\exp\left\{2\int_{1-\lambda_1^{-\varepsilon/2}}^1 -\partial_1V_1^1(0,s)\, ds\right\}.
\end{equation}
We also denote
\begin{equation}\notag
K_1(t)=e^{\int_t^1-\partial_1\bar V_1^1(0,s)\, ds}.
\end{equation}
It follows from \eqref{dx-V11} that
\begin{equation}\notag
\partial_1 V_1^1(0,t)=\partial_1 \bar V_1^1(0,t)=
-g_1(t)\lambda_1^{\gamma+2+\frac12(1-t)\lambda^{\varepsilon/2}}\tilde \phi_0^1 \tilde \phi_0^2\approx g_1(t)\lambda_1^{\gamma+2}
\end{equation}
where we used $\frac12\leq |\tilde \phi_0^1|\leq 2$ and \eqref{phi2-size} again. Moreover, for large enough $\lambda>0$, we have $\lambda_1\gg \lambda$ and hence $1-\lambda^{-\frac12\varepsilon}+\lambda_1^{-\varepsilon}< 1-\lambda_1^{-\frac12\varepsilon}$. Thus the definition of $g_1$ in \eqref{g1} indicates $g_1(t)=1$ for $t\in [1-\lambda_1^{-\frac12\varepsilon},1]$.
As a consequence, it holds
\begin{equation}\label{lam2}
\lambda_2\approx \exp\left\{2\int_{1-\lambda_1^{-\varepsilon/2}}^1 \lambda_1^{\gamma+2}\, ds\right\} \approx e^{2\lambda_1^{\gamma+2-\frac12\varepsilon}}.
\end{equation}
We observe $\lambda_2\gg \lambda_1$ for $\gamma=-2+\varepsilon$. On the other hand, we have
\begin{equation}\label{k1}
K_1(t)\approx e^{\lambda_1^{2+\gamma}\int_t^1 g_1(s)\, ds}\approx e^{(1-t)\lambda_1^{2+\gamma}}, \ \ 1-\lambda_1^{-\frac12\varepsilon}\leq t\leq 1.
\end{equation}
It is clear to see $K_1(t)\leq \lambda_2^{\frac12}$.

We are ready to construct $W_2$. First let $\bar W_2$ be the solution to the linear equation with an ending time data 
\begin{equation}\label{eq-bar-b2}
\begin{split}
\partial_t\bar W_2-\bar V_1\cdot\nabla \bar W_2=-\bar W_2\cdot\nabla \bar V_1\\
\bar W_{2,\mbox{end}}(\vec x):=\bar W_2(\vec x,1)=(0,0,\bar W_2^3(\vec x,1))
\end{split}
\end{equation}
with
\begin{equation}\label{bar-b13}
\bar W_2^3(\vec x,1)=-\lambda_2^{\gamma}\sin(\lambda_2 x)f(\lambda_2^{\gamma_1}x)\sin(\lambda_2 y)f(\lambda_2^{\gamma_2}y).
\end{equation}
Let $\phi_1(\vec x, t_0,t)$ be the solution to the flow map \eqref{flow} with $\bar J_0$ replaced by $\bar V_1$, that is,
\begin{equation}\label{flow1}
\begin{cases}
\partial_t\phi_1^i(\vec x, t_0, t)=-\bar V_1^i(\phi_1(\vec x, t_0, t), t), \ \ i=1,2,3\\
\phi_1(\vec x, t_0, t_0)= \vec x.
\end{cases}
\end{equation} 

Based on the properties of $\bar V_1$ and the particular ending data of \eqref{eq-bar-b2}, we notice again that the first equation of \eqref{eq-bar-b2} is equivalent to 
\begin{equation}\notag
\partial_t\bar W_2^i-\bar V_1\cdot\nabla \bar W_2^i=-\bar W_2^i\partial_i \bar V_1^i, \ \ i=1,2,3.
\end{equation}
That is, each component $\bar W_2^i$ satisfies a linear equation. Therefore, applying Lemma \ref{le-solution} to \eqref{eq-bar-b2}, \eqref{bar-b13} and \eqref{flow1}, we have the solution of \eqref{eq-bar-b2} 
\begin{equation}\notag
\bar W_2^i(\vec x,t)=\bar W^i_{2,\mbox{end}}(\phi_1(\vec x,t,1))\exp\left\{\int_1^t-\partial_i\bar V_1^i\left(\phi_1(\vec x, t,s),s\right)\, ds \right\}
\end{equation}
with
\begin{equation}\label{supp2}
\supp\left( \bar W_2(\vec x,t)\right)\subset\left\{\vec x=(x,y,z): |x|\leq 2\lambda_2^{-\gamma_1}K_1(t), \ |y|\leq 2\lambda_2^{-\gamma_2}K_1(t) \right\}.
\end{equation}
Therefore,
\[\bar W_2^1(\vec x,t)\equiv 0, \ \ \bar W_2^2(\vec x,t)\equiv 0, \ \ \bar W_2^3(\vec x,t)=\bar W_{2,\mbox{end}}^3(\phi_1(\vec x,t,1)).\]

Define the smooth temporal cutoff function $g_2$ 
\begin{equation}\label{g1}
g_2(t)=
\begin{cases}
0, \ \ \ \ t\leq 1-\lambda_1^{-\varepsilon/2}\\
(t-1+\lambda_1^{-\varepsilon/2})\lambda_{2}^{\varepsilon}, \ \ \ 1-\lambda_1^{-\varepsilon/2}\leq t\leq 1-\lambda_1^{-\varepsilon/2}+\lambda_2^{-\varepsilon}\\
1, \ \ \ \ 1-\lambda_1^{-\varepsilon/2}+\lambda_2^{-\varepsilon}\leq t\leq 1
\end{cases}
\end{equation}
and the final perturbation 
\[W_2(\vec x, t)=g_2(t)\bar W_2(\vec x,t).\]
The new solution $B_2=B_1+W_2$ with the current density $J_2=\nabla\times B_2$ satisfies the forced equation
\begin{equation}\notag
\partial_tB_2-J_2\cdot \nabla B_2=-B_2\cdot\nabla J_2+F_1+F_2
\end{equation}
where, by using \eqref{iter-b1} and \eqref{eq-bar-b2} 
\begin{equation}\notag
\begin{split}
F_2=&\ g_2'(t)\bar W_2+g_2(t)\left[(\bar V_1-J_1)\cdot\nabla \bar W_2-\bar W_2\cdot\nabla (\bar V_1-J_1) \right]\\
&+g_2(t)\left[-(\nabla\times \bar W_2)\cdot\nabla B_1+B_1\cdot\nabla(\nabla\times \bar W_2) \right]\\
&+g_2^2(t)\left[-(\nabla\times \bar W_2)\cdot\nabla \bar W_2+\bar W_2\cdot\nabla(\nabla\times \bar W_2) \right].
\end{split}
\end{equation}
Noticing 
\[(\nabla\times \bar W_2)\cdot\nabla \bar W_2=0, \ \ \bar W_2\cdot\nabla(\nabla\times \bar W_2)=0,\]
and 
\[ B_1\cdot\nabla(\nabla\times \bar W_2)=B_1^3\partial_z(\nabla\times \bar W_2)=0, \] 
\[ \bar W_2\cdot\nabla(\bar V_1-J_1)=\bar W_2^3\partial_z (\bar V_1-J_1)=0,\]
the force $F_2$ is essentially
\begin{equation}\label{F2}
F_2=g_2'(t)\bar W_2+g_2(t)(\bar V_1-J_1)\cdot\nabla \bar W_2-g_2(t)(\nabla\times \bar W_2)\cdot\nabla B_1.
\end{equation}

\medskip

We can verify the flow map $\phi_1$ satisfies similar estimates as in Lemma \ref{le-phi}.
\begin{Lemma}\label{le-phi1}
Let $\gamma_2>\frac12$. Denote 
\begin{equation}\notag
\phi_1^1(\vec x,t,1)=xK_1(t) \tilde \phi_1^1(y,t).
\end{equation}
For $t\in[1-\lambda_1^{-\frac12\varepsilon}, 1]$ and $|y|\leq 2\lambda_2^{-\gamma_2}K_1(t)$, we have
\begin{equation}\notag
\begin{split}
\frac12\leq |\tilde \phi_1^1(y,t)|\leq 2, \ \ |\partial_2\tilde \phi_1^1(y,t)|\leq 8 \lambda^{4+\gamma}\lambda_2^{-\gamma_2}K_1(t), \ \
|\partial_2^2\tilde \phi_1^1(y,t)|\leq 8\lambda_1^{4+\gamma},\\
\frac{1}{2K_1(t)}\leq |\partial_2\phi_1^2(\vec x,t,1)|\leq\frac{2}{K_1(t)},\ \ 
|\partial_2^2\phi_1^2(\vec x,t,1)|\leq \lambda_1^{4+\gamma}\lambda_2^{-\gamma_2}K_1(t).
\end{split}
\end{equation}
\end{Lemma}
Moreover, the following H\"older estimates can be justified analogously.
\begin{Lemma}\label{le-b2-holder}
Let $k\in \mathbb N$ and $\beta\in(0,1)$. We have
\[ \|\bar W_2(t)\|_{C^{k,\beta}}\lesssim \lambda_2^{\gamma+k+\beta}, \ \ \mbox{for} \ 1-\lambda_2^{-\frac12\varepsilon}\leq t\leq 1.\]
\end{Lemma}
The proof is similar to that of Lemma \ref{le-b1-holder} and thus omitted.
\begin{Lemma}\label{le-F2-holder}
We have for $1-\lambda_2^{-\frac12\varepsilon}\leq t\leq 1$
\[ \|F_2(t)\|_{C^\alpha}\lesssim \lambda_2^{1+\gamma+\alpha}.\]
\end{Lemma}
\pf
Since
\[g_2'(t)\sim \lambda_2^{\varepsilon}, \ \ \|\bar W_2\|_{C^\alpha}\lesssim \lambda_2^{\gamma+\alpha}, \ \ \|\bar W_1\|_{C^\alpha}\lesssim \lambda_1^{\gamma+\alpha}\]
and $B_1=B_0+g_1(t)\bar W_1$, it follows
\[|g_2'(t)\bar W_2|\lesssim \lambda_2^{\varepsilon+\gamma}, \ \ |g_2(t)(\nabla\times \bar W_2)\cdot\nabla B_1| \lesssim \lambda_1^{\gamma+1}\lambda_2^{\gamma+1}.\]
On the other hand, it has
\[\bar V_1-J_1=(\bar V_1-V_1)+(V_1-J_1)=(\bar V_1-V_1)-J_0.\]
It is easy to see from \eqref{V11} that $V_1(0,t)=0$. Therefore, invoking \eqref{bar-V11} we deduce
\begin{equation}\notag
\bar V_1-V_1=\frac12 x^2\partial_x^2 V_1((0,y,0),t)+ \mbox{lower order terms},
\end{equation}
which implies
\begin{equation}\notag
|\bar V_1-V_1|\lesssim  |x|^2\|\bar W_1\|_{C^3}\lesssim  |x|^2 \lambda_1^{\gamma+3}.
\end{equation}
It then follows from \eqref{supp2}, $|x|\leq 2\lambda_2^{-\gamma_1}K_1(t)\lesssim 1$ and $|y|\leq 2\lambda_2^{-\gamma_2}K_1(t)\lesssim 1$ for $\gamma_1, \gamma_2>\frac12$, and consequently
\begin{equation}\notag
\begin{split}
|g_2(t)(\bar V_1-J_1)\cdot\nabla \bar W_2|&\lesssim \left(|\bar V_1-V_1|+|J_0|\right)|\nabla \bar W_2|\\
&\lesssim (\lambda_1^{\gamma+3}|x|^2+\lambda^{\gamma+3}|x||y|)\lambda_2^{1+\gamma}\\
&\lesssim \lambda_1^{\gamma+3}\lambda_2^{1+\gamma}.
\end{split}
\end{equation}
In the end we claim from \eqref{F2}
\begin{equation}\notag
\begin{split}
\|F_2(t)\|_{C^{\alpha}}\lesssim&\ \|\bar W_2\|_{C^\alpha}+\|(\bar V_1-J_1)\cdot\nabla \bar W_2\|_{C^\alpha}+\|(\nabla\times \bar W_2)\cdot\nabla B_1\|_{C^\alpha}\\
\lesssim&\ \lambda_2^{\varepsilon+\gamma+\alpha}+\lambda_2^{1+\gamma+\alpha}\\
\lesssim&\ \lambda_2^{1+\gamma+\alpha}.
\end{split}
\end{equation}

\cbdu

\medskip

\subsection{Iterative construction of perturbation $W_{q+1}$ for $q\geq 2$}
\label{sec-bq}

We have demonstrated in Subsection \ref{sec-b2} that, employing the magnetic field profile with zero horizontal components which generates a hyperbolic current flow near origin, the iterating process can be carried on. In this subsection, we state the iteration for the $q$-th step, $q\geq 2$.

For any $q\geq 1$, let $B_q=(0,0, B_q^3)$ be a solution to the forced equation
\begin{equation}\label{eq-bq}
\partial_tB_q-J_q\cdot \nabla B_q=-B_q\cdot\nabla J_q+\sum_{j=1}^qF_j, \ \ \ J_q=\nabla\times B_q.
\end{equation}
Assume $B_q^3$ does not dependent on the $z$ variable. We further assume for $1-\lambda_2^{-\frac12\varepsilon}\leq t\leq 1$,
\begin{equation}\label{assu}
 \|W_q(t)\|_{C^\alpha}\lesssim \lambda_q^{\gamma+\alpha},\ \ \
 \|F_q(t)\|_{C^\alpha}\lesssim \lambda_q^{1+\gamma+\alpha}
 \end{equation}
with $W_q=B_q-B_{q-1}$. Denote $V_q=\nabla\times W_q$.
To construct the next perturbation $W_{q+1}$, we use the modified current density
$\bar V_q=(\bar V_q^1, \bar V_q^2, \bar V_q^3)$ defined by
\begin{equation}\label{bar-Vq}
\begin{cases}
\bar V_q^1=x\partial_x V_q^1((0,y,0),t),\\
\bar V_q^2= V_q^2((0,y,0),t),\\
\bar V_q^3=0.
\end{cases}
\end{equation}

Define iteratively 
\begin{equation}\notag
\lambda_{q+1}=\exp\left\{2\int_{1-\lambda_q^{-\varepsilon/2}}^1 -\partial_1V_q^1(0,s)\, ds\right\},
\end{equation}
\begin{equation}\notag
K_q(t)=e^{\int_t^1-\partial_1\bar V_q^1(0,s)\, ds},
\end{equation}
and the smooth temporal cutoff function $g_{q+1}$ 
\begin{equation}\label{gq}
g_{q+1}(t)=
\begin{cases}
0, \ \ \ \ t\leq 1-\lambda_q^{-\varepsilon/2}\\
(t-1+\lambda_q^{-\varepsilon/2})\lambda_{q+1}^{\varepsilon}, \ \ \ 1-\lambda_q^{-\varepsilon/2}\leq t\leq 1-\lambda_q^{-\varepsilon/2}+\lambda_{q+1}^{-\varepsilon}\\
1, \ \ \ \ 1-\lambda_q^{-\varepsilon/2}+\lambda_{q+1}^{-\varepsilon}\leq t\leq 1
\end{cases}
\end{equation}

Let $\bar W_{q+1}$ be the solution to the linear equation with an ending time data 
\begin{equation}\label{eq-bar-bq}
\begin{split}
\partial_t\bar W_{q+1}-\bar V_q\cdot\nabla \bar W_{q+1}=-\bar W_{q+1}\cdot\nabla \bar V_q\\
\bar W_{q+1,\mbox{end}}(\vec x):=\bar W_{q+1}(\vec x,1)=(0,0,\bar W_{q+1}^3(\vec x,1))
\end{split}
\end{equation}
with
\begin{equation}\label{bar-bq3}
\bar W_{q+1}^3(\vec x,1)=-\lambda_{q+1}^{\gamma}\sin(\lambda_{q+1} x)f(\lambda_{q+1}^{\gamma_1}x)\sin(\lambda_{q+1} y)f(\lambda_{q+1}^{\gamma_2}y).
\end{equation}
Let $\phi_q(\vec x, t_0,t)$ be the solution to the flow map 
\begin{equation}\label{flowq}
\begin{cases}
\partial_t\phi_q^i(\vec x, t_0, t)=-\bar V_q^i(\phi_q(\vec x, t_0, t), t), \ \ i=1,2,3\\
\phi_q(\vec x, t_0, t_0)= \vec x.
\end{cases}
\end{equation} 

The particular form of $\bar V_q$ in \eqref{bar-Vq} and the ending data of \eqref{eq-bar-bq} ensure that the first equation of \eqref{eq-bar-bq} is equivalent to 
\begin{equation}\notag
\partial_t\bar W_{q+1}^i-\bar V_q\cdot\nabla \bar W_{q+1}^i=-\bar W_{q+1}^i\partial_i \bar V_q^i, \ \ i=1,2,3.
\end{equation}
Consequently we have the explicit formula of the solution.
\begin{Lemma}\label{le-solution-q}
The solution of \eqref{eq-bar-bq}-\eqref{bar-bq3} is given by
\begin{equation}\label{pert-q}
\bar W_{q+1}^i(\vec x, t)=\bar W_{q+1,\mbox{end}}^i\left( \phi_q(\vec x, t,1)\right) \exp\left\{ \int_1^t-\partial_i \bar V_q^i\left(\phi_q(\vec x, t, s),s \right) \, ds\right\}
\end{equation}
for $i=1,2,3$.
Moreover, we have
\begin{equation}\notag
\supp\left( \bar W_{q+1}(\vec x,t)\right)\subset\left\{\vec x=(x,y,z): |y|\leq 2\lambda_{q+1}^{-\gamma_2}K_q(t) \right\}.
\end{equation}
\end{Lemma}
The flow map $\phi_q$ defined by \eqref{flowq} has the following properties.
\begin{Lemma}\label{le-phiq}
Let $\gamma_2>\frac12$. Denote 
\begin{equation}\notag
\phi_q^1(\vec x,t,1)=xK_q(t) \tilde \phi_q^1(y,t).
\end{equation}
For $t\in[1-\lambda_{q}^{-\frac12\varepsilon}, 1]$ and $|y|\leq 2\lambda_{q+1}^{-\gamma_2}K_q(t)$, we have
\begin{equation}\notag
\begin{split}
\frac12\leq |\tilde \phi_q^1(y,t)|\leq 2, \ \ |\partial_2\tilde \phi_q^1(y,t)|\leq 8 \lambda^{4+\gamma}\lambda_{q+1}^{-\gamma_2}K_q(t), \ \
|\partial_2^2\tilde \phi_q^1(y,t)|\leq 8\lambda_q^{4+\gamma},\\
\frac{1}{2K_q(t)}\leq |\partial_2\phi_q^2(\vec x,t,1)|\leq\frac{2}{K_q(t)},\ \ 
|\partial_2^2\phi_q^2(\vec x,t,1)|\leq \lambda_q^{4+\gamma}\lambda_{q+1}^{-\gamma_2}K_1(t).
\end{split}
\end{equation}
\end{Lemma}
\begin{Lemma}\label{le-bq-holder}
For $1-\lambda_{q+1}^{-\frac12\varepsilon}\leq t\leq 1$ and $\alpha\geq 0$, the H\"older estimate
\[ \|\bar W_{q+1}(t)\|_{C^{\alpha}}\lesssim \lambda_{q+1}^{\gamma+\alpha}\]
holds.
\end{Lemma}
Lemmas \ref{le-solution-q}, \ref{le-phiq} and \ref{le-bq-holder} can be proved analogously as for Lemmas \ref{le-solution}, \ref{le-phi} and \ref{le-b1-holder} respectively, with minor modifications. 

Since $\bar W_{q+1, \mbox{end}}^1=\bar W_{q+1, \mbox{end}}^2=0$ and $\partial_3\bar V_q^3=0$, the solution formula \eqref{pert-q} implies
\[\bar W_{q+1}(\vec x,t)=\left(0,0, \bar W_{q+1,\mbox{end}}^3\left( \phi_q(\vec x, t,1)\right)  \right)\]
which does not depend on the $z$ variable.
Define the new perturbation $W_{q+1}(\vec x,t)=g_{q+1}(t)\bar W_{q+1}(\vec x,t)$ and solution 
\[B_{q+1}=B_q+W_{q+1}(\vec x,t)=\left(0,0, B_q^3+g_{q+1}(t)\bar W_{q+1}^3(\vec x,t) \right)\]
which also does not depend on the $z$ variable. Denote $J_{q+1}=\nabla\times B_{q+1}$.
Combining \eqref{eq-bq} and \eqref{eq-bar-bq} yields the forced equation of $B_{q+1}$,
\begin{equation}\label{eq-bq1}
\partial_tB_{q+1}-J_{q+1}\cdot \nabla B_{q+1}=-B_{q+1}\cdot\nabla J_{q+1}+\sum_{j=1}^{q+1}F_j
\end{equation}
with 
\begin{equation}\notag
\begin{split}
F_{q+1}=&\ g_{q+1}'(t)\bar W_{q+1}+g_{q+1}(t)\left[(\bar V_q-J_q)\cdot\nabla \bar W_{q+1}-\bar W_{q+1}\cdot\nabla (\bar V_q-J_q) \right]\\
&+g_{q+1}(t)\left[-(\nabla\times \bar W_{q+1})\cdot\nabla B_q+B_q\cdot\nabla(\nabla\times \bar W_{q+1}) \right]\\
&+g_{q+1}^2(t)\left[-(\nabla\times \bar W_{q+1})\cdot\nabla \bar W_{q+1}+\bar W_{q+1}\cdot\nabla(\nabla\times \bar W_{q+1}) \right].
\end{split}
\end{equation}
Again, since the only non-zero component (vertical) of $\bar W_{q+1}$ does not depend on $z$ variable, it has
\[(\nabla\times \bar W_{q+1})\cdot\nabla \bar W_{q+1}=0, \ \ \bar W_{q+1}\cdot\nabla(\nabla\times \bar W_{q+1})=0,\]
\[ B_q\cdot\nabla(\nabla\times \bar W_{q+1})=B_q^3\partial_z(\nabla\times \bar W_{q+1})=0, \] 
\[ \bar W_{q+1}\cdot\nabla(\bar V_q-J_q)=\bar W_{q+1}^3\partial_z (\bar V_q-J_q)=0.\]
Hence the force $F_{q+1}$ is reduced to
\begin{equation}\label{Fq}
F_{q+1}=g_{q+1}'(t)\bar W_{q+1}+g_{q+1}(t)(\bar V_q-J_q)\cdot\nabla \bar W_{q+1}-g_{q+1}(t)(\nabla\times \bar W_{q+1})\cdot\nabla B_q.
\end{equation}

Applying the assumption \eqref{assu} and Lemma \ref{le-bq-holder}, we can show
\begin{Lemma}\label{le-Fq-holder}
For $1-\lambda_{q+1}^{-\frac12\varepsilon}\leq t\leq 1$ and $\alpha\geq 0$, the estimate
\[ \|F_{q+1}(t)\|_{C^\alpha}\lesssim \lambda_{q+1}^{1+\gamma+\alpha}\]
holds.
\end{Lemma}

\bigskip

\section{Proof of Theorem \ref{thm}}
\label{sec-proof}
Let $\lambda>0$ be a sufficiently large constant. For any small constant $\varepsilon>0$, fix the parameters $\gamma=-2+\frac12\varepsilon$ and $\gamma_1=\gamma_2=\frac12+\varepsilon$.
With such choice, starting from the profile $B_0=(0,0, B_0^3)$ with $B_0^3$ defined in \eqref{b0-3}, we can obtain iteratively a sequence of solutions $B_q$ satisfying the forced electron MHD \eqref{eq-bq} for $q\geq 1$ in the previous section. Denote the limiting solution by
\[B(\vec x, t)=\lim_{q\to\infty}B_q(\vec x, t), \ \ \ J(\vec x, t)=\nabla\times B(\vec x, t), \] and the external forcing \[F(\vec x, t)=\lim_{q\to\infty}\sum_{j=1}^q F_j(\vec x, t).\] 
We claim that the limits exist on the time interval $[0,1)$. Indeed, for any $0\leq t\leq 1-\lambda_q^{-\frac12\varepsilon}$, the choice of the temporal cutoff function in \eqref{gq} ensures that
\[B_p(\vec x, t)=B_q(\vec x, t), \ \ \sum_{j=1}^p F_j(\vec x, t)=\sum_{j=1}^q F_j(\vec x, t), \ \ \mbox{for any} \ \ p\geq q,\]
and hence 
\[B(\vec x, t)=B_q(\vec x, t), \ \ \lim_{p\to\infty}\sum_{j=1}^p F_j(\vec x, t)=\sum_{j=1}^q F_j(\vec x, t)\]
on $0\leq t\leq 1-\lambda_q^{-\frac12\varepsilon}$. Thus the limit magnetic field $B(t)$ is a classical solution to the forced equation \eqref{emhd} on the time interval $[0,1)$.
Moreover, it follows from Lemma \ref{le-bq-holder} and Lemma \ref{le-Fq-holder} that for $0\leq t\leq 1$,
\begin{equation}\notag
\begin{split}
\|B(t)\|_{C^{2-\varepsilon}}\leq \|B_0\|_{C^{2-\varepsilon}}+\lim_{q\to\infty}\sum_{j=1}^q \|\bar W_j(t)\|_{C^{2-\varepsilon}}\lesssim \lambda^{\gamma+2-\varepsilon}+\sum_{j=1}^\infty \lambda_j^{\gamma+2-\varepsilon}\lesssim 1,\\
\|F(t)\|_{C^{1-\varepsilon}}\leq \lim_{q\to\infty}\sum_{j=1}^q \|F_j(t)\|_{C^{1-\varepsilon}}\lesssim \sum_{j=1}^\infty \lambda_j^{1+\gamma+1-\varepsilon}\lesssim 1
\end{split}
\end{equation}
Since $B_q(t)$ and $F_q(t)$ are compactly supported on a small spatial set, it is trivial to see $B(t), F(t)\in L^2(\mathbb R^3)$. Therefore item (i) in Theorem \ref{thm} is verified. 

To show item (ii), we take a sequence of time $t_q=1-\lambda_{q}^{-\frac12\varepsilon}$ which converges to $t=1$ as $q\to\infty$. At $t_q$ we have 
\[\nabla\times J(t_q)=-\Delta B(t_q)=-\Delta B_0-\sum_{j=1}^{q}g_j(t)\Delta \bar W_j(t_q)\]
where we used the temporal cutoff property again. Observing that $\Delta \bar W_{q}$ is the leading order term and 
\[\left|g_{q}(t)\Delta \bar W_{q}(t_q)\right|\approx \lambda_q^{\varepsilon} \lambda_q^{\gamma+2}\approx \lambda_q^{\frac32\varepsilon}\]
in view of \eqref{gq}, \eqref{bar-bq3} and \eqref{pert-q}. Thus, as $q\to \infty$, $t_q\to 1$ and $ \lambda_q^{\frac32\varepsilon}\to\infty$. It completes the proof of (ii).

In the end, we prove (iii). It follows from the analysis above
\begin{equation}\notag
\begin{split}
\int_0^1 \|\nabla\times J(t)\|_{L^\infty}\, dt\approx& \sum_{j=1}^\infty \|\nabla\times J(t_j)\|_{L^\infty}(t_{j+1}-t_j)\\
\approx& \sum_{j=1}^\infty \lambda_j^{\frac32\varepsilon}\left(\lambda_j^{-\frac12\varepsilon}-\lambda_{j+1}^{-\frac12\varepsilon}\right)\\
\approx& \sum_{j=1}^\infty \lambda_j^{\varepsilon}.
\end{split}
\end{equation}
Therefore we have $\int_0^1 \|\nabla\times J(t)\|_{L^\infty}\, dt=\infty$.

\bigskip

\section{Construction of shear profile and proof of Theorem \ref{thm2}}
\label{sec-thm2}

In this section we provide another scenario of blowup for the forced electron MHD relying on shear current profile construction and prove Theorem \ref{thm2}. 
 
Let $f$ be the smooth cutoff function defined in Subsection \ref{sec-profile}. We then take the compactly supported magnetic field $B_0=(0,0, B_0^3)$ with
\begin{equation}\notag
\bar B_0^3=-\lambda_0^\gamma \sin(\lambda_0 x)f(\lambda_0^{\frac12}x)
\end{equation}
for a large constant $\lambda_0>0$ and parameter $\gamma=-2+\frac12\varepsilon$. It is easy to see $\div B_0=0$. 
The current density produced by $B_0$ is 
\[\bar J_0=\nabla\times \bar B_0=(0, -\partial_x \bar B_0^3, 0)\]
with 
\begin{equation}\notag
\bar J_0^1=\bar J_0^3=0, \ \
\bar J_0^2=-\partial_x \bar B_0^3= \lambda^{1+\gamma}\cos(\lambda x)f(\lambda^{\gamma_1}x)
+\lambda^{\gamma+\gamma_1}\sin(\lambda x)f'(\lambda^{\gamma_1}x).
\end{equation}
It is clear that 
\[\bar J_0\cdot\nabla \bar B_0=0, \ \ \ \bar B_0\cdot\nabla \bar J_0=0.\]
Thus $\bar B_0$ is a stationary solution to the electron MHD \eqref{emhd} without forcing.

For a fixed constant $b>1$, define the frequency number 
\[\lambda_q=\lambda_0^{b^q}, \ \ q\geq 1.\] 
For large enough $\lambda_0>0$, we note $\lambda_{q+1}\gg \lambda_q$. For each $q\geq 1$, we define the sheer magnetic field
\begin{equation}\label{b-sheer}
\bar B_q=(0,0, \bar B_q^3) \ \ \mbox{with} \ \ \bar B_q^3=-\lambda_q^\gamma \sin(\lambda_q x)f(\lambda_q^{\frac12}x).
\end{equation}
Its current density is 
\[\bar J_q=\nabla\times \bar B_q=(0, -\partial_x \bar B_q^3, 0)\]
with 
\begin{equation}\label{j-sheer}
\begin{split}
\bar J_q^1&=\bar J_q^3=0, \\
\bar J_q^2&=-\partial_x \bar B_q^3= \lambda_q^{1+\gamma}\cos(\lambda_q x)f(\lambda_q^{\gamma_1}x)
+\lambda_q^{\gamma+\gamma_1}\sin(\lambda_q x)f'(\lambda_q^{\gamma_1}x).
\end{split}
\end{equation}
Let $\bar W_{q+1}$ be the solution to the linear equation with ending data specified at $t=1$
\begin{equation}\label{sheer-bq}
\begin{split}
\partial_t\bar W_{q+1}-\bar J_q\cdot\nabla \bar W_{q+1}=-\bar W_{q+1}\cdot\nabla \bar J_q\\
\bar W_{q+1}(\vec x,1)=(0,0,\bar W_{q+1}^3(\vec x,1))
\end{split}
\end{equation}
with
\begin{equation}\label{sheer-bq3}
\bar W_{q+1}^3(\vec x,1)=-\lambda_{q+1}^{\gamma}\sin(\lambda_{q+1} x)f(\lambda_{q+1}^{\frac12}x).
\end{equation}
First we observe, due to the form of $\bar J_q$ in \eqref{j-sheer},
\[\bar W_{q+1}\cdot\nabla \bar J_q=(0, \bar W_{q+1}^1\partial_x \bar J_q^2, 0).\]
Hence the equation in \eqref{sheer-bq} is equivalent to the system
\begin{equation}\label{sheer-bq2}
\begin{split}
\partial_t\bar W_{q+1}^1-\bar J_q\cdot\nabla \bar W_{q+1}^1&=0,\\
\partial_t\bar W_{q+1}^2-\bar J_q\cdot\nabla \bar W_{q+1}^2&=-\bar W_{q+1}^1\partial_x \bar J_q^2,\\
\partial_t\bar W_{q+1}^3-\bar J_q\cdot\nabla \bar W_{q+1}^3&=0.
\end{split}
\end{equation}
Since $\bar W_{q+1}^1(\vec x,1)=0$, the first (linear) equation of \eqref{sheer-bq2} implies $\bar W_{q+1}^1(\vec x,t)\equiv 0$. Thus the second equation of \eqref{sheer-bq2} is also purely transport, i.e.
\[\partial_t\bar W_{q+1}^2-\bar J_q\cdot\nabla \bar W_{q+1}^2=0\]
which together with $\bar W_{q+1}^2(\vec x,1)=0$ indicates $\bar W_{q+1}^2(\vec x,t)\equiv 0$.

In the end, we note 
\[\bar J_q\cdot\nabla \bar W_{q+1}^3=\bar J_q^1\partial_x \bar W_{q+1}^3+\bar J_q^2\partial_y \bar W_{q+1}^3+\bar J_q^3\partial_z \bar W_{q+1}^3=\bar J_q^2\partial_y \bar W_{q+1}^3.\]
Thus the third equation of \eqref{sheer-bq2} reduces to 
\[\partial_t\bar W_{q+1}^3-\bar J_q^2\partial_y \bar W_{q+1}^3=0\]
whose solution satisfying the ending data \eqref{sheer-bq3} is 
\begin{equation}\notag
\bar W_{q+1}^3(\vec x,t)=-\lambda_{q+1}^{\gamma}\sin(\lambda_{q+1} x)f(\lambda_{q+1}^{\frac12}x).
\end{equation}

Summarizing the analysis above gives the solution to \eqref{sheer-bq}
\begin{equation}\label{wq}
\bar W_{q+1}(\vec x,t)=\left(0,0,  -\lambda_{q+1}^{\gamma}\sin(\lambda_{q+1} x)f(\lambda_{q+1}^{\frac12}x)\right).
\end{equation}

Let $g_q(t)$ be the smooth temporal cutoff functions constructed in \eqref{gq}. Define the perturbations
\[W_{q+1}(\vec x,t)=g_{q+1}(t) \bar W_{q+1}(\vec x,t)\]
and solutions iteratively
\[B_{q+1}(\vec x,t)=B_q(\vec x,t)+W_{q+1}(\vec x,t)=B_0(\vec x)+\sum_{j=1}^{q+1}g_{j}(t) \bar W_{j}(\vec x,t).\]
We can verify that 
\[\nabla\times ((\nabla\times B_{q+1})\times B_{q+1})=B_{q+1}\cdot\nabla J_{q+1}-J_{q+1}\cdot \nabla B_{q+1} =0.\]
Hence $B_{q+1}$ satisfies the forced electron MHD 
\begin{equation}\label{bfq}
\partial_tB_{q+1}+\nabla\times ((\nabla\times B_{q+1})\times B_{q+1})=\sum_{j=1}^{q+1}F_j
\end{equation}
with $F_{j}= g_{j}'(t)\bar W_{j}$.

It follows from \eqref{wq} and \eqref{gq} that for $0\leq t\leq 1-\lambda_q^{-\frac12\varepsilon}$
\[B_p(\vec x, t)=B_q(\vec x, t), \ \ \sum_{j=1}^p F_j(\vec x, t)=\sum_{j=1}^q F_j(\vec x, t), \ \ \mbox{for any} \ \ p\geq q;\]
and for $1-\lambda_q^{-\frac12\varepsilon}\leq t\leq 1$,
\[\|W_{q+1}(t)\|_{C^\alpha}\approx \lambda_{q+1}^{\gamma+\alpha}, \ \ \|F_{q+1}(t)\|_{C^\alpha}\approx \lambda_{q+1}^{\frac12\varepsilon+\gamma+\alpha}.\]
Similar argument as in Section \ref{sec-proof} shows that the limiting magnetic field $B(\vec x, t)=\lim_{q\to\infty} B_q(\vec x, t)$ is a solution to the forced electron MHD \eqref{emhd} on the time interval $[0,1)$ with the external forcing
\[F(\vec x,t)=\sum_{j=1}^{\infty}F_j.\]
Moreover, since $\gamma=-2+\frac12\varepsilon$, 
\begin{equation}\notag
\begin{split}
\|B(t)\|_{C^{2-\varepsilon}}\leq \|B_0\|_{C^{2-\varepsilon}}+\sum_{j=1}^\infty \|\bar W_j(t)\|_{C^{2-\varepsilon}}\lesssim \lambda_0^{-2+\frac12\varepsilon+2-\varepsilon}+\sum_{j=1}^\infty \lambda_j^{-2+\frac12\varepsilon+2-\varepsilon}\lesssim 1,\\
\|F(t)\|_{C^{2-\frac32\varepsilon}}\leq \sum_{j=1}^\infty \|F_j(t)\|_{C^{2-\frac32\varepsilon}}\lesssim \sum_{j=1}^\infty \lambda_j^{\frac12\varepsilon-2+\frac12\varepsilon+2-\frac32\varepsilon}\lesssim 1.
\end{split}
\end{equation}
Thus $B(t)\in C^{2-\varepsilon}(\mathbb R^3)\cap L^2(\mathbb R^3)$ and $F(t)\in C^{2-\frac32\varepsilon}(\mathbb R^3)\cap L^2(\mathbb R^3)$ on $[0,1]$, where the $L^2$ boundedness is guaranteed by the compact support property.  The items (ii), (iii) and (iv) can be obtained analogously as in Section \ref{sec-proof}.



\end{document}